\documentclass[10pt]{elsart3-1} 
\usepackage{amssymb}
\usepackage[english,francais]{babel}
\usepackage[figuresright]{rotating}
\usepackage{amsmath, amsfonts}
\newcommand \RR    {\mathbb{R}} 
\newcommand \EE    {\mathbb{E}} 

\newcommand \Scal    {\mathcal S} 
\newcommand \Lcal    {\mathcal{L}} 
 
\newcommand \del   {\partial} 

\newtheorem{e-proposition}[theorem]{Proposition} 

\newtheorem{e-definition}[theorem]{Definition\rm}

\renewcommand{\theequation}{\arabic{equation}}
\def\og{\leavevmode\raise.3ex\hbox{$\scriptscriptstyle\langle\!\langle$~}}
\def\fg{\leavevmode\raise.3ex\hbox{~$\!\scriptscriptstyle\,\rangle\!\rangle$}}

\newcommand \vs {\vskip.3cm}

\newcommand \be {\begin{equation}}
\newcommand \ee {\end{equation}}
\newcommand \bei {\begin{itemize}}
\newcommand \eei {\end{itemize}}
    
\journal{the Acad\'emie des sciences}
\begin{document} 
\centerline{}
\begin{frontmatter}
 
\selectlanguage{english}
\title{An algorithm (CoDeFi) for 
overcoming the curse of dimensionality in mathematical finance}

\selectlanguage{english}
\author[authorlabel1]{Philippe G. LeFloch}
\ead{contact@philippelefloch.org}
\and 
\author[authorlabel2]{Jean-Marc Mercier}
\ead{jean-marc.mercier@mpg-partners.com}

\address[authorlabel1]{Laboratoire Jacques-Louis Lions \& Centre National de la Recherche Scientifique 
\\
Universit\'e Pierre et Marie Curie, 4 Place Jussieu, 75252 Paris, France. } 

\address[authorlabel2]{MPG-Partners, 136 Boulevard Haussmann, 75008 Paris, France.}

\medskip

\begin{abstract}
\selectlanguage{english}  
We present an algorithm (CoDeFi) which overcomes the curse of dimensionality (CoD) in scientific computations and, especially, in mathematical finance (Fi). Our method applies a broad class of partial differential equations such as Kolmogorov-type equations and, for instance, the Black and Scholes equation. As a main feature, our algorithm allows one to solve partial differential equarions in large dimensions and provides a general framework for stochastic problems. In insurance or finance applications, the number of dimensions corresponds to the number of risk sources and it is crucial to have a numerical method that remains robust and reliable in large dimensions.




\vskip.15cm 

\selectlanguage{francais}

\noindent{\bf R\'esum\'e.} 
Nous pr\'esentons un algorithme (CoDeFi) pour r\'esoudre le probl\`eme de la mal\'ediction dimensionnelle (CoDe, en anglais) pour le calcul scientifique et, en particulier, pour aborder des probl\`emes provenant la finance (Fi). 
Notre m\'ethode s'applique \`a une large classe d'\'equations aux d\'eriv\'ees partielles telles que les \'equations de Kolmogorov et, par exemple, l'\'equation de Black et Scholes. La sp\'ecificit\'e de notre m\'ethode est de permettre la r\'esolution d'\'equations aux d\'eriv\'ees partielles en grandes dimensions et de fournir un environnement g\'en\'eral pour les probl\`emes stochastiques. Dans les applications dans le secteur de l'assurance ou de la finance, le nombre de dimensions correspond au nombre de sources de risques et il est crucial d'avoir une m\'ethode num\'erique qui reste robuste et pr\'ecise en grandes dimensions.
\end{abstract}
\end{frontmatter}

\selectlanguage{francais}

\section*{Version fran\c{c}aise abr\'eg\'ee} 
 
\vskip-.3cm  La \textit{mal\'ediction dimensionnelle} (CoD, ce terme provenant d'un article de 1957 de R.E. Bellman) appara\^\i t lorsque le temps de calcul augmente  exponentiellement avec le nombre de dimensions, qui correspond au nombre de sources de risque (ou sous-jacents) dans les applications en math\'ematiques financi\`eres. Nous proposons dans cette Note une approche g\'en\'erale \`a ce probl\`eme ouvert depuis plusieurs d\'ecennies. Notre approche est bas\'ee sur une combinaison de techniques provenant de la th\'eorie des \'equations aux d\'eriv\'ees partielles: des trajectoires de Monte-Carlo, une m\'ethode num\'erique classique qui n'est pas "maudite" (c'est-\`a-dire non affect\'ee par CoD), et des maillages de calcul mobiles. De tels maillages permettent de r\'esoudre les \'equations de Kolmogorov en utilisant des m\'ethodes non structur\'ees conjointement avec une technique de transport optimal. Nous utilisons \'egalement une calibration ``z\'ero-erreur'' qui nous permet d'assurer une description pr\'ecise de la dynamique du processus sous-jacent. Toutes ces techniques sont regroup\'ees dans un environnement de calcul scientifique, que nous avons d\'evelopp\'e (CoDeFi) et qui peut \^etre vu comme un environnement de mesure de risque tr\`es g\'en\'eral, bas\'e l'algorithme pr\'esent\'e en d\'etail dans \cite{PLF-JMM-1}-\cite{PLF-JMM-2}.

Cette note pr\'esente ces technologies et propose un premier benchmark pour le cas multi-dimensionnel, en le compl\'etant jusqu\`a la 64 \`eme dimension (ou 64 sources de risque). A notre connaissance, peu d'autres technologies pourraient compl\'eter ce test. La quantification optimale \cite{BaPa} ou l'analyse par ondelette \cite{MaNiSc} pourraient \^etre utilis\'es jusqu\''\`a disons $10$ dimensions. Au del\`a, les techniques type Monte-Carlo am\'ericain \cite{LS} pourraient peut-\^etre fournir des bornes inf\'erieures, ces m\'ethodes \'etant connues pour fournir des exercices sub-optimaux. 
Nous concluons en notant que ce m\^eme environnement de calcul est utilis\'e pour des simulations dans le cadre de probl\`emes hyperboliques nonlin\'eaires \cite{PLF-JMM-1}.


\vskip-.65cm\selectlanguage{english} 

\section*{English version} 

\section{Introduction and main strategy}
 
\vskip-.3cm   {\it 1.1 Fokker-Planck-Kolmogorov problems.} 
We denote by $\mu=\mu(t,x)$ ($t \ge 0$, $x \in \RR^D$) the probability density measure of a stochastic process. 
Here $D$ denotes the number of dimensions, i.e. usually the number of risk sources, also called underlyings, in the 
applications in mathematical finance. We suppose that this process defines a Markov chaining process and follows a Fokker-Planck equation of the form 
\be \label{FP}
	\del_t \mu + \Lcal \mu = 0, 
\ee 
where $\Lcal(t,\cdot)$ is a (usually unknown) positive differential operator of parabolic type. In fact, we are interested in the adjoint of Fokker-Planck equations, that is, backward Kolmogorov equations with unknown $P=P(t,x) \in \RR^M$ (defined for $0\le t \le T$ and $x \in \RR^D$): 
\be \label{KE}
	\del_t P - \Lcal^* P = 0, \quad P(T,\cdot) = P_T \in \RR^M, 
\ee 
where $\Lcal^*$ denotes the adjoint and the equation is submitted to a terminal Cauchy data $P_T=P_T(x)$ at some (future) time $T > 0$. This is typically the case of flows (e.g. cash or asset based flows) in mathematical finance. We also consider here nonlinear systems extending \eqref{KE}, which we can formulate as an  optimal-stopping problem
\be \label{KEO}
	\text{either } \del_t P - \Lcal^* P=0 \quad \text{or else both} \quad \del_t P - \Lcal^* P \ge 0 \text{ and } \Scal(t,x,P) = 0. 
\ee
Here, the mapping $\Scal=\Scal(t,x,P): \RR^{D+1+M} \to \RR^P$ is interpreted as a ``strategy''.  In this context, the following two rather distinct problems arise:

\textbf{1- The calibration problem}. This consists in finding out a probability measure $\mu(t,\cdot)$, satisfying a list of integral constraints at some (future) times denoted by $T_1, T_2, \ldots, T_I >0$. These constraints have the form 
\be \label{Co}
	\int_{\RR^D} P^i( \cdot) \, \mu (T_i, \cdot) = C_i, \qquad i = 1, \ldots, I,
\ee
in which the functions $P^i$ and the functions $C_i$ (refered to as observables) are prescribed. Starting from a local volatility model "\`a la Dupire", several methods are available in order to handle this (undetermined) problem; we refer for instance \cite{Ho} for a review of commonly used methods. 

\textbf{2- The valuation problem}. The probability measure $\mu(t,\cdot)$ being now given, the second problem of interest consists in expliciting $\mathcal{L}$ and solving the equation \eqref{KEO}, leading to the output measure computed from the surface $\{t,x,P(t,x)\}$. There exists also several numerical methods devoted to solving \eqref{KE} and \eqref{KEO}; see for instance \cite{BT} for a survey of Monte-Carlo-type methods (that are able to handle linear problems like \eqref{KE}) and PDE (partial differential equation) lattice-based methods (able to handle also the nonlinear problems \eqref{KEO}).

Note that we need to compute the whole ``surface'' $\{t,x,P(t,x)\}$ to compute properly the nonlinear solution to \eqref{KEO}. Even for linear problems like \eqref{KE}, this surface contains all the information required for risk measurement applications
(of operational or regulatory nature). 
To our knowledge, only a PDE approach is capable to compute the surface $\{t,x,P(t,x)\}$. However, 
this approach cannot be used for problems that can involve hundreds of dimensions, say.
No method yet is capable to tackle the calibration and valuation problems above in a systematic and general manner. For instance, optimal quantization \cite{BaPa} or wavelet analysis  \cite{MaNiSc} could be used, but for only up to about $10$ dimensions. Above this threshold, American-Monte Carlo methods \cite{LS} could be used and provide lower bounds to solutions of \eqref{KEO}, but these methods are known to compute sub-optimal exercising. The lack of numerical methods able to handle large dimensions was identified first by Bellman in 1957 and, since then is refered to as  the \textit{Curse of dimensionality}.  
In other words, a general, swiss-knife method is important for applications
and, 
if the same methodology could be used by practitionners, then risk measures could be compared in a meaningful way confidently aggregated in order to deduce macro-economical indicators.

{\it 1.2 Application to mathematical finance.} For instance, one could treat the scalar strategy $\Scal(t,x,P) = P - (x-K)^+$ (with $Q^+:= \max(Q,0)$) corresponding to an optimal strategy for American-type options having strike $K \in \RR$. The vector-valued function $P=P(t,x)$ 
in general determines the \textit{fair} value of contracts under the strategy $\Scal$, when the underlying stochastic process is worth $x \in \RR^D$ at a future time $t \ge 0$.
The method we now propose provides a very general framework since all the stochastic models usually defined by prescribing specific stochastic processes, (such as normal, log-normal, stochastic volatility, local volatility or local correlation models...) fit into the Fokker-Planck setting \eqref{FP}. Moreover, the applications in mathematical finance involve contraints of the form \eqref{Co} which we can handle here. As stated above, all classical risk measurements can be deduced from the knowledge of the solutions of \eqref{KEO}.  For instance for operational measures, $P^t:= \int_{\RR^D} P(t,\cdot) \mu(t,\cdot)$ corresponds to the \textit{fair} value if $t=0$ (and future value if $t >0$) of contracts viewed from today ($t=0$), since $\int_{\RR^D} \nabla P(t,\cdot) \mu(t,\cdot)$ computes its hedge, where $\nabla:= \big(\frac{\partial}{\partial x_d}\big)_{d=1..D}$ is the gradient operator. Similarly for regulatory type measures, a VaR (Value At Risk) can be expressed as a quantile $\int_{ \{x: P(t, x) < \alpha P^t\} } P(t,\cdot) \mu(t,\cdot)$, where $t$ is the VaR time horizon, since $0 < \alpha < 1$ is a confidence threshold. A CVA (Credit Value Adjustment) can be expressed from the knowledge of the function of $t \mapsto \int_{\RR^D } P^+(t,\cdot) \mu(t,\cdot)$.

{\it 1.3 Purposes and main ideas}. The method CoDeFi which we propose in this Note is based on several novel numerical techniques, which lead us to a framework for handling the calibration problem as well as the valuation problem, in a high dimensional setting.
Our CoDeFi algorithm is based on two main ingredients. Our first idea is a \textbf{localisation principle}: we use a change of variable to localize the system of equations \eqref{KEO} into the unit cube, denoted $\Lambda = [0,1]^D$. To properly introduce this change of variable, we recall a standard result from the theory of optimal transport (cf.~Villani \cite{Vil}), going back to Brenier \cite{Brenier}) allowing one to regard the quantile of a probability measure as a change of variable:
\be \label{quantile}
 \mu(t,\cdot) = S(t,\cdot)_\# m, \qquad S(t,\cdot) = \nabla h:\Lambda \mapsto \RR^D, \quad h \text{ convex}. 
\ee
where $S(t,\cdot)_\# m$  stands for the pull-back of the Lebesgue measure $m$ on $\Lambda$. For our purpose, we think of $S(t,\cdot)$ as a map inducing a natural change of variable within Monte-Carlo methods. Hence, $S(t,\cdot)$ can be used  together with a random generator and allows us to sample some underlying process at any time $t$, by writing 
\be \label{MC}
	\EE^t(P, \mu ):= \int_{\RR^D} P(\cdot) \mu(t, \cdot) \simeq \frac{1}{N} \sum_{n=1, \ldots, N} P\circ S(t,Y_n), \quad n=1, \ldots, N,
\ee
where $P$ is any function, and the matrix $Y = \{ Y_n \in \Lambda \}_{n=1, \ldots, N} \in \RR^{N\times D}$ consists of {\it random}   vectors. Our localization principle then consists in considering the following \textit{transported} (into the unit cube) version of \eqref{KEO}:
\be \label{KEOT}
	\text{either } \big( \del_t P - \Lcal^* P \big) \circ S=0 
 \quad \text{or else both } \quad \big( \del_t P - \Lcal^* P \big) \circ S \ge 0 \text{ and } \Scal(t, S,P \circ S) = 0
\ee
One important incentive for using the change of variables above is to transform the Kolmogorov problem \eqref{KE} into a \textit{self-adjoint} problem. (See \cite{PLF-JMM-2} for the details.)

Our second ingredient is an \textbf{unstructured mesh technique} in order to solve the problem \eqref{KEOT} in a robust and inexpensive way. Our motivation here is clear: random sampling methods (such as Monte-Carlo ones \eqref{MC}) are not sensible to the Curse of Dimensionality. Hence, if a PDE-based method is able to use the sampling vectors $Y$ as a mesh, then the resulting solver will not be affected by the CoD. In the rest of this note, we build upon these two main ideas and provide   details on the implementation of the CoDeFi algorithm.

 \section{The calibration and valuation algorithms}

\vskip-.3cm  
{\it 2.1 Generation of the computational grid}. 
There are two essentials parameters, central in all notations, and main limiting factors in term of computational cost, that is, the number $\textit{D}$ of underlyings or risk sources, and the number \textit{N} of grid points.
The overall complexity of the CoDeFi algorithm is bounded by $\mathcal{O}((N+D)^3)$ operations. For such a number of operations, the relative error commited to solve Kolmogorov equations \eqref{KEO} is $\epsilon = \mathcal{O}(1/N) \%$. To have some figures in mind, CoDeFi solves today Kolmogorov equations \eqref{KEO} with $D=64$, $N=512$, with time of order 15 minutes, for a precision of order $1/512 \sim 0.2 \%$, using one core of the eight of an Intel 4770 processor.

\begin{figure}[h]
\centering
\begin{tabular}{ccc}

\begin{minipage}{1.5in}
\includegraphics[height=1.6in,width=1.9in]
{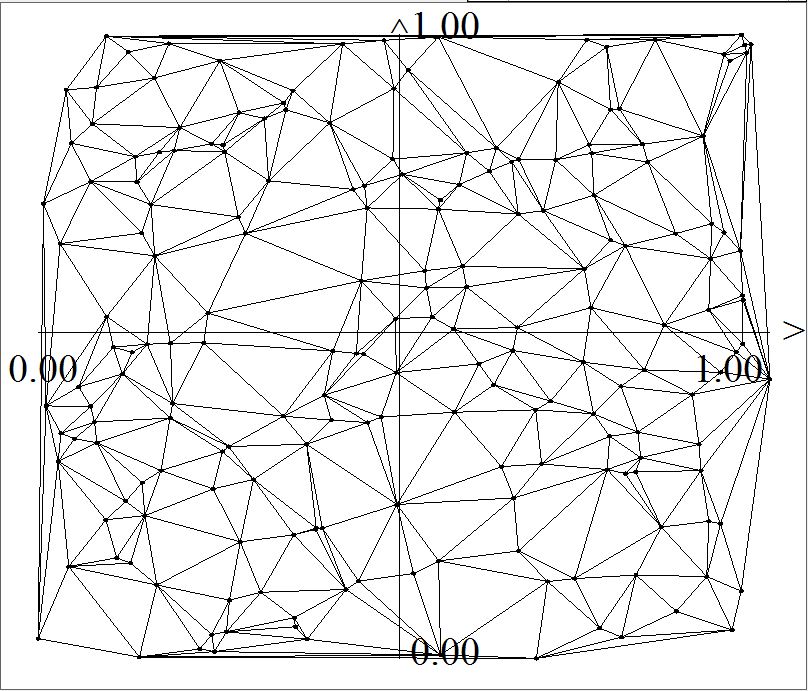}
\centerline{Quasi Random Sequence}
\end{minipage}

\hskip2.01cm

\begin{minipage}{1.5in}
\includegraphics[height=1.6in,width=1.9in]
{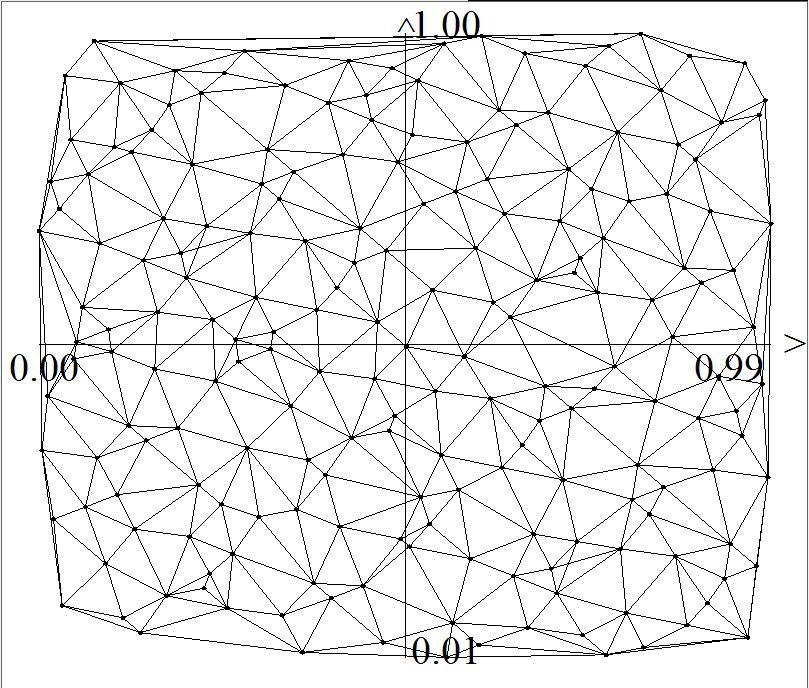}
\centerline{ Sobol sequence}
\end{minipage}
\hskip2.01cm

\begin{minipage}{1.5in}
\includegraphics[height=1.6in,width=1.9in]
{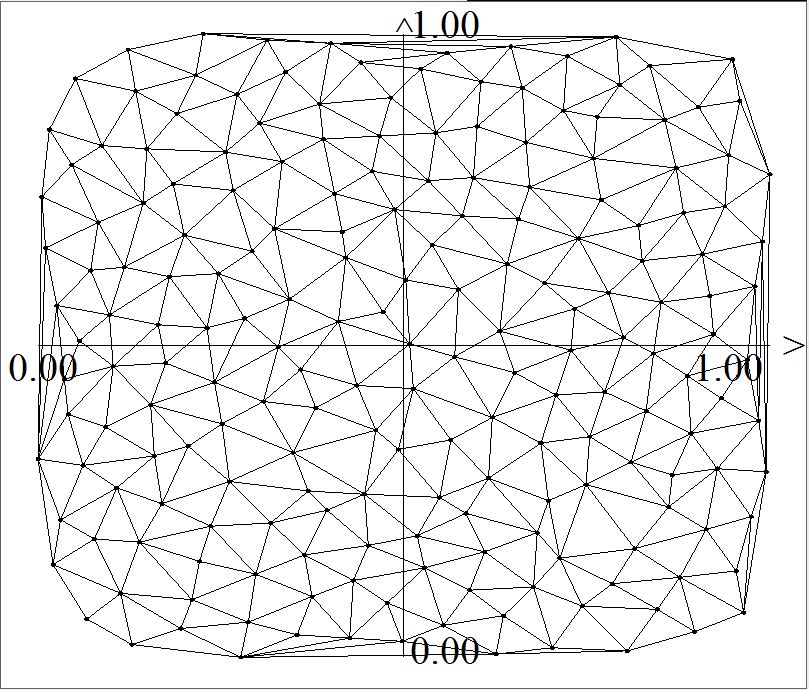}
\centerline{Optimal Discrepancy Sequence}
\end{minipage}

\end{tabular}
\caption{Different grids with Delaunay meshes} 
\label{BEE-Sobol}
\end{figure}

We specify here the choice of the $N$ sampling-set $Y = \{ Y_n \in \Lambda \}_{n=1, \ldots, N} \in \RR^{N\times D}$, used as a mesh for the PDE engine. CoDeFi is primarily a Monte-Carlo method, and can use any kind of sampling-set of the uniform law over $\Lambda$ as grid. For instance Figure \ref{BEE-Sobol} shows the plot of a two-dimensional grid of 200 points, with Delaunay mesh, using the Mersenne Twister generator mt19937, see \cite{MN}. We recall that Pseudo-random generator owns a statistical convergence rate of order $\epsilon = \mathcal{O}(1/\sqrt{N}) \%$. This error holds for Monte-Carlo integration formula of kind \eqref{MC} and bounded variations functions $P$.

Another alternative is to use quasi-random sequences. Figure \ref{BEE-Sobol} also shows the plot of a 200 points grid using a Sobol, see \cite{Sobol}, low-discrepancy generator. We recall that Sobol generators have a  convergence rate of order $\epsilon = \mathcal{O}((\ln N)^D /{N }) \%$, for bounded variations functions $P$ and formula \eqref{MC}.

Finally, since mesh repartition is an important feature for PDE methods, we propose to use special sequences, that we call optimal discrepancy sequences, to generate the grid.  Such sequences reaches the above cited rate of convergence, of order $\epsilon = \mathcal{O}(1 / N) \%$, for finite sums of convex and concave functions $P$ and formula \eqref{MC}. However such sequences require $\mathcal{O}(N^2D)$ operations to compute. Since PDE methods works usually with few grid points ($N=512$ is enough), we prefer sacrifying computational time to mesh quality.

\vs

{\it 2.2 The calibration algorithm}. 
We are now given a mesh as a matrix $Y = \{ Y_n \in \Lambda \}_{n=1, \ldots, N} \in \RR^{N\times D}$, determined as in the previous section.
As quoted above, the calibration algorithm consists in finding a probability measure $\mu(t,\cdot)$ that fits constraints \eqref{Co}. Due to Brenier result, this is equivalent to finding out a quantile $S(t,\cdot)$ - see \eqref{quantile} - fitting constraints.  Let $\mu_0(t,\cdot)$ a given, smooth "prior" probability measure - in practice, we use normal or log-normal processes as priors - and denotes $S_0(t,\cdot)$ its quantile, known explictely. From a numerical point of view, we solve the calibration problem finding a matrix  $S(t) = \{ S^n(t):= \nabla h(t,Y_n) \}_{n=1, \ldots, N} \in \RR^{N\times D}$, $h$ convex, satisfying for each distinct time $T_i$,
$$	
\inf_{S(T_i) \in \RR^{ND}} \|S(T_i) -  S_0(T_i)\|_2,\quad \text{const. } \frac{1}{N} \sum_{n=1, \ldots, N} P^i(T_i, S^n(T_i))  = C^i, \quad i=1, \ldots, I,
$$
where $\| \cdot\|_2$ denotes the standard Frobenius norm of matrixes and $S_0(t) = \{ S_0(t,Y_n) \}_{n=1, \ldots, N} \in \RR^{N D}$. The solution of this problem can be computed with classical Lagrangian methods, and is not sensitive to the CoD. We can also show that this problem is equivalent to finding a discrete measure minimizing $\inf_{ \mu(T_i):= \frac{1}{N} \sum\delta_{S_n(T_i)} } \| \mu(T_i) -  \mu_0(T_i,\cdot)\|_{W^2}$ under constraints \eqref{Co}, where $  \| \cdot \|_{W^2}$ denotes the \textit{Wasserstein} distance and $\delta_X$ the Dirac mass weighting a point $X \in \RR^D$. Starting from distinct calibrated times $T_i$, standard bootstrap arguments and interpolation allows to retrieve the whole surface $t \mapsto S(t) \in \RR^{ND}$ for any times if needed.

\begin{table}
\caption{Market values of European Call SX5E Mat 3M. Spot 3064.03}
\label{SX5E}
\begin{tabular*}{\columnwidth}{@{}l@{\extracolsep{\fill}}r@{\extracolsep{\fill}}r@{\extracolsep{\fill}}r@{\extracolsep{\fill}}r@{\extracolsep{\fill}}r@{\extracolsep{\fill}}r@{\extracolsep{\fill}}r@{\extracolsep{\fill}}l@{\extracolsep{\fill}}c}
& & \multicolumn{6}{c}{}
 \\ [-6pt]
{Strike $\alpha_i$ }  & 0.8 & 0.9 & 0.95 & 0.975 & 1 & 1.025 & 1.05 & 1.1 & 1.2 \\
\hline
{Value $C^i$ } & 559.2  & 292.6 & 180.7  & 133.6  & 93.76  & 61.59 & 37.99 & 11.34 & 0.31 \\
\hline
\end{tabular*}
\end{table}

\begin{figure}
 \centerline{\includegraphics[width=1.5in]{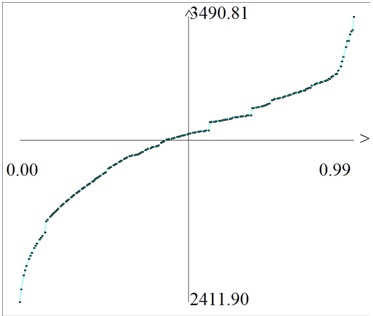}}
\caption{Calibrated SX5E quantile for $N=256$, $D=1$.}
\label{CSX5E}
\end{figure}
Note that the constraints setting \eqref{Co} is a quite general one. In particular, we can use a set of constraints to match exactly classical statistical measures. For instance
\begin{itemize}
  \item \textit{Expectations}: we can match any of the $D$ marginal expectations $T \mapsto F_d(T):= \int x_d \mu(T,\cdot)$ of any distributions using the constraint $F_d(T,x) = x_d$ in \eqref{Co};
  \item \textit{Variances}: we can match any of the $D$ marginal variances $T \mapsto V_d(T)$ of any distributions using the constraints $V_d(T,x) = \Big( x_d  - F_d(T) \Big)^2 $ in \eqref{Co};
  \item \textit{Correlations}: we can match any correlation matrix $T \mapsto C_{d_1,d_2}(T), d_1,d_2 = 1..D,$ using the constraints $C_{d_1,d_2}(T,x) = \frac{( x_{d_1}  - F_{d_1}(T)) ( x_{d_1}  - F_{d_2}(T))}{V_{d_1}(T) V_{d_2}(T)}$, $1 \le d_1 < d_2 \le D$ in \eqref{Co};
\end{itemize}

{\it 2.3 Application to mathematical finance.} 
For instance, consider Table \ref{SX5E}, describing real constraints, corresponding to quotes of call options written over index $SX5E$ having maturity $3$ months (3M), that is $P^i(T_i,x) = (x-\alpha_i S(0))^+$, $T_ i = 3M$, where $\alpha_i$ is a percentage of today value $S(0)$. Then Figure \ref{CSX5E} shows a quantile calibrated to these quotes at time 3M. 
Constraints like \eqref{Co} cannot be always matched. Even in an operational context, \textit{arbitrage} occurs in market data and an important feature of our approach is that, since we solve a constrained problem, the calibration remains stable but cannot match constraints. Indeed, we can use this calibration algorithm to detect market arbitrages within a large variety of situations, and might also propose an explanation to highlight them.

{\it 2.4 An example.} 
To illustrate this feature, as well as the capacity of this calibration algorithm to handle higher dimensions and numerous constraint efficiently, consider a independant log-normal process $X_d=e^{Z_d /10 }$, $d=1.D$, where $Z_d$ denotes a standard normal process.  Table \eqref{EBO10Y} presents, for $D=1,4,16,64$, $N=32,128,512$, the computation of the expectation \eqref{MC}, for the function $P(T,x) = ( |x|_\infty - K)^+$, with $K=1$, (called a Best-of option), where $|x|_\infty = \sup_{d=1..D}(|x_d|)$, and $T=10 Y$. Table \ref{EBO10Y}-line MC- presents reference computations using $N=1048576 \times D$ sequence of a pseudo-random Mersenne twister MT19937, that are confident with a relative error estimated to $0,1 \% \sim 1/ \sqrt{1048576}$. The others columns presents the same computations, but using calibrated sequences $S(T)$ as described in this section, matching all expectations, all variances, and the whole correlation identity matrix, that is 2144 constraints for $D=64$. Indeed, we noticed that this calibration procedure accelerates the convergence rate of Monte-Carlo sampling \eqref{MC}, as could be expected.

\begin{table}
\caption{Price of European Best-of options MAT 10Y}
\label{EBO10Y}
\begin{tabular*}{\columnwidth}{@{}l@{\extracolsep{\fill}}r@{\extracolsep{\fill}}l@{\extracolsep{\fill}}c@{\extracolsep{\fill}}c}
& & \multicolumn{2}{c}{}\\ [1pt]
 \\ [-6pt]
{ } &{D=1} &  {D=4} &  {D=16}&  {D=64} \\
\hline
{ MC N = 1048576} & 0.12578 & 0.35941   & 0.68796  & 1.0166\\
\hline
{ N=32} & 0.128275  & 0.340435  & 0.678092 & 0.927398\\
{ N=128} & 0.126521  & 0.349573  & 0.693397 & 0.982611\\
{ N=512} & 0.125921 & 0.359632  & 0.688982  & 1.0144\\
\hline
\end{tabular*}
\end{table}

{\it 2.5 The valuation algorithm}. We now  present some numerical results for \eqref{KEO}, once the calibration step has been performed. In this valuation phase, we suppose that the quantile $S(t,\cdot)$ of the probability measure $\mu(t,\cdot)$ is known, and the sampling set $S(t) \in \RR^{N \times D}$ is computed as described in the previous section.

Consider the linear equation \eqref{KE} and denotes $P(t) \sim P(t,S(t)) \in \RR^{N \times M}$ an approximated solution, where $t$ belongs to a fixed time-grid. Let $s \le t$ two consecutive times of this time grid. Then \eqref{KE} is approximated using a classical Crank-Nicolson approach, and the solution is computed accordingly to
\be \label{KED}
	P(s) =  \Pi^{ (t,s)}P(t),  \quad \Pi^{ (t,s)}:= \Big( \pi^{ (t,s)}_{n,m} \Big)_{n,m=1, \ldots, N} \in \RR^{N \times N},
\ee
where the matrix $\Pi^{ (t,s)}$ (the approximated generator to the Kolmogorov equation \eqref{KE}) is computed explicitly. 

This matrix can be interpreted in a Markov-chaining process setting: $\pi^{ (t,s)}_{n,m}$ is the probability that the stochastic process jumps from the state $S^t_n \in \RR^D$ to the state $S^s_m \in \RR^D$. Hence an important property of this matrix is to be \textit{bi-stochastic}, reflecting the fact that the underlying process defines a \textit{martingale} process. Obviously, this bi-stochastic property is obtained after projecting  the quantile $t\mapsto S(t,\cdot)$ into a proper space of functions. In particular, consider any constraints $P^i(T_i, \cdot)$ to the underlying process, see \eqref{Co}. Then, this \textit{bi-stochastic} property implies the following discrete conservation of expectations for any time 
\be \label{Cal}
	C_i = \frac{1}{N} \sum_{n=1, \ldots, N} P^i(t,S^n(t)), \quad 0 \le t \le T_i
\ee
For instance, considering the example treated in Table \ref{EBO10Y}, \eqref{Cal} implies that the correlation matrix of the underlying process perfectly fit the identity matrix for any time $0 \le t \le T$. Solutions to the nonlinear problem \eqref{KEO} are approximated adding nonlinear terms to the linear solution \eqref{KED} at each time, as follows: 
$
	P(s) = \max \big( \Pi^{ (t,s)}P(t),  \Scal(s, S(s), \Pi^{ (t,s)}P(t)) \big)$.
Hence, once the transition matrices $\Pi^{ (t,s)}$ computed, computing \eqref{KEO} is mainly a matter of matrix-matrix multiplication. Indeed, for vector valued solutions of size $M = ND$ (i.e. $\sim$ 32000 for $N= 512$, $D=64$), no significant slow down are noticed. Moreover, matrix-matrix multiplication is usually a quite optimized code, and this backward step can be easily distributed.

\section{Numerical experiments}

\vskip-.5cm

\begin{figure}[h]
\centering
\begin{tabular}{ccc}

\begin{minipage}{1.5in}
\includegraphics[height=1.6in,width=1.9in]
{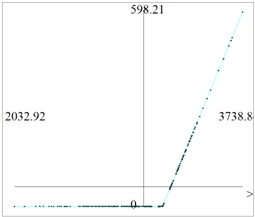}
\centerline{Payoff $T = 3M$}
\label{AMCALL3M}
\end{minipage}

\hskip2.01cm

\begin{minipage}{1.5in}
\includegraphics[height=1.6in,width=1.9in]
{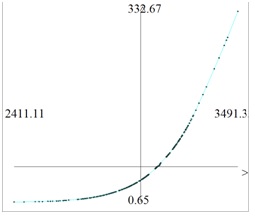}
\centerline{ $T = 1M$}
\label{AMCALL1M}
\end{minipage}
\hskip2.01cm

\begin{minipage}{1.5in}
\includegraphics[height=1.6in,width=1.9in]
{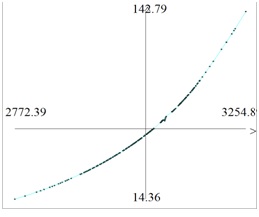}
\centerline{$T = 6D$}
\label{AMCALL6D}
\end{minipage}

\end{tabular}
\caption{$P(T,S)$ at different retropopagation time, $N=256$} 
\end{figure}

{\it 3.1 One-dimensional experiments.}
Our first test is a one dimensional test, that consists in solving \eqref{KEO}, with the same data as presented in Table \ref{SX5E}, but using the strategy $\Scal(t,x,P) = P - (x-K_i)^+$, corresponding to an american option exercise. We considered the data set in Table \ref{SX5E} to calibrate the underlying process, hence the quantile is the same as plotted in Figure \ref{CSX5E}. The qualitative properties of the solution are plotted in a serie of three figures, illustrating the retropagation steps for the special case $K_i = 3064$. The first one, Figure \ref{AMCALL3M}, represents the initial conditions at time $T=3M$. It plots more precisely the surface $\{S^{3M}_n, P(3M,S^{3M}_n) \}_{n=1..256}$. The second one, figure \eqref{AMCALL1M}, plots the solution at time $t=1M$, $\{S^{1M}_n, P(1M,S^{1M}_n) \}_{n=1..256}$. The third one  plots the solution at time $t=6D$, $\{S^{6D}_n, P(6D,S^{6D}_n) \}_{n=1..256}$, that is the last computation time. 

\begin{table}
\caption{American Call SX5E Mat 3M. Spot 3064.03}
\label{SX5EAM}
\begin{tabular*}{\columnwidth}{@{}l@{\extracolsep{\fill}}r@{\extracolsep{\fill}}r@{\extracolsep{\fill}}r@{\extracolsep{\fill}}r@{\extracolsep{\fill}}l@{\extracolsep{\fill}}c}
& & \multicolumn{4}{c}{}\\ [1pt]
 \\ [-6pt]
{Strike \%} & {\it Call values} & {\it N=16} & {\it N=64} & {\it N=256}& {\it N=1024} \\
\hline
0.8 & 559.224 & 650.54  & 650.54 & 650.54 & 650.54 \\
0.9 & 292.61 & 346.98  & 347.17 & 347.17 & 347.20 \\
0.95 & 180.78  & 209.90  & 210.27 & 210.29 & 210.32 \\
0.975 & 133.59  & 151.00  & 152.27 & 152.22 & 152.23 \\
1 & 93.76  & 103.16  &104.19 & 104.33 & 104.37 \\
1.025 & 61.59 & 66.53  & 67.19 & 67.37 & 67.40 \\
1.05 & 37.99 & 40.23  & 40.68 & 40.79 & 40.82 \\
1.1 & 11.34 & 19.47  & 11.83 & 11.92 & 11.91 \\
1.2 & 0.31 & 0.31  &0.31 & 0.31 &0.32 \\
\hline
comp. time &  & 0.02 s.   & 0.06 s. & 1.52 s. & 77 s. \\
\end{tabular*}
\end{table}
Expectations of the solution at time $T=6D$, that corresponds to fair american options prices, are presented in Table \ref{SX5EAM}, with differents number of points $N$ - we recall that $N$ drives the precision of the computations - and different strikes $K_i$, but the same as used for the calibration process. It has already been noted that these numerical schemes are very stable and accurate \cite{JMM1}, as confirmed by this table. Computation time is also confirming the $N^3$ behavior of the algorithm. For $D=1$, we could specialize the algorithm to get a complexity of order $\mathcal{O}(N)$ (as in \cite{JMM1}), but we prefer to use the same algorithm, independently of the dimension $D$.

{\it 3.2 Multi-dimensional experiments}.
Illustrating the multi-dimensional case is not an easy task, since there is no available benchmark. Hence, we propose here a  benchmark, for which the error can be estimated. Consider the European Best-of options in Table \ref{EBO10Y}, and consider the strategy $\Scal(t,x,P) = P - (|x |_\infty-1)^+$ for times $t=1Y,\ldots, 10Y$, corresponding to a virtual product with annual exercise (Bermudean style option). Since the initial condition $P(T,x) = (|x |_\infty-1)^+$ is convex, the solution remains convex, and the strategy should never be exercised, i.e. $P(t,x) \ge (|x |_\infty-1)^+$, where $P$ is solution to the linear problem \eqref{KE}. The expectations of the solution (i.e. the fair price) can be estimated using results in Table \ref{EBO10Y}.

Table \ref{ECO10Y} presents the expection of the solution $\frac{1}{N} \sum_{n=1, \ldots, N} P^i(t,S^n(t))$, for $t= 6D$, computed with our algorithm. The error is estimated as a percentage with respect to \eqref{EBO10Y}, as well as an indicative computational time. These results show that the strategy is slightly exercized, i.e. the solution to the linear problem \eqref{KE} does not satisfy $P(t,x) \ge (|x |_\infty-1)^+$ as expected. Indeed, this indicates that, even if we are using a numerical scheme that is convergent and stable, the solution presents small oscillations which could probably by suppressed by improving this step of the method. 


\begin{table}
\caption{Price of Bermudean Best-of options MAT 10Y}
\label{ECO10Y}
\begin{tabular*}{\columnwidth}{@{}l@{\extracolsep{\fill}}c@{\extracolsep{\fill}}c@{\extracolsep{\fill}}c@{\extracolsep{\fill}}c}
& & \multicolumn{2}{c}{}
 \\ [-6pt]
{ } &{D=1} &  {D=4} &  {D=16}&  {D=64} \\
\hline
{ MC N = 1048576} & 0.12578 & 0.35941   & 0.68796  & 1.0166\\
\hline
{ N=32} & 0.128275 \quad (0 \%, 0.14 s.) & 0.40047 \quad (8.1 \%, 0.21s.) & 0.684418 \quad (0.46 \%,2.3s.) & 0.927398 \quad (0.36 \%, 152s.)\\
{ N=128} & 0.126521 \quad (0 \%, 1.30 s.) & 0.38120 \quad (4.3 \%, 2.49s.) &0.711784 \quad (1.30 \%, 14.4s.) & 0.985387 \quad (0.14 \%, 304s.)\\
{ N=512} & 0.125921 \quad (0 \%, 121 s.) & 0.39757  \quad (5.0 \%, 148s.) & 0.699329  \quad (0.74 \%, 280s.) & 1.01721\quad (0.13 \%, 1665s.)\\
\hline
\end{tabular*}
\end{table}

{\it 3.3 Concluding  remarks.} In this note we presented a methodology to solve parabolic / hyperbolic equations set in high dimensions. Even if it remains perfectible, the proposed can already confidently be used in an industrial environment to compute risk measurements. There are numerous potential applications linked to the curse of dimensionality. We identified some of them in the insurance industry, but mainly in the financial sector. Above pricing and hedging, allowing to issue new financial products, more adapted to clients needs, or market abitraging, we are adressing risk measurements: indeed, the test realized in this paper shown that CoDeFi could compute accurate Basel-based regulatory measures as RWA (Risk Weighted Asset) or CVA (Credit Valuation Adjustment) on a single laptop, whereas farm of thousand of computers are used on a daily basis today with approximation methods. Finally, such a framework could be helpful in systemic risk measurements, that can be accurately modeled with high dimensional Kolmogorov equations like \eqref{KEO}.


\end{document}